\documentclass[11pt]{article}
\newcommand{ \newsection}[1]{ \setcounter{equation}{0} \section{ #1} }

\def\keywords{ \if@twocolumn
\section*{Keywords}
\else \small
\begin{center}
{ \bf Keywords\vspace{-.5em}\vspace{0pt}}
\end{center}
\center \fi}
\def\endkeywords{ \if@twocolumn\else\endcenter\fi}
\usepackage{epsfig}
\usepackage{graphicx}
\begin{document}
\title{\bf Solution of Seventh Order Boundary Value Problems using Adomian Decomposition
Method}
%%%%%%%%%%%%%%%%%%%%%%%%%%%%%%%%%%%%%%%%%%%%%%%%%%%%%%%%%%%%%%%%%%%%%%%%%%
\author{Shahid S. Siddiqi \thanks{
Department of Mathematics, University of the Punjab, Lahore 54590,
 Pakistan.   ~~~~~~~~~~~~~~~~~~~~~~~~~~~~~~~~~~~~~~~~~~~~~~~~~~~~~~~~Email: shahidsiddiqiprof@yahoo.co.uk},
%  Ghazala Akram \thanks{
%Department of Mathematics, University of the Punjab, Lahore 54590,
% Pakistan.   ~~~~~~~~~~~~~~~~~~~~~~~~~~~~~~~~~~~~~~~~~~~~~~~~~~~~~~~~Email: toghazala2003@yahoo.com},
 ~Muzammal Iftikhar \thanks{ Department of Mathematics, University of the Punjab, Lahore 54590,
 Pakistan.   ~~~~~~~~~~~~~~~~~~~~~~~~~~~~~~~~~~~~~~~~~~~~~~~~~~~~~~~~Email: miftikhar@hotmail.com}}
\date{}
\maketitle
%%%%%%%%%%%%%%%%%%%%%%%%%%%%%%%%%%%%%%%%%%%%%%%%%%%%%%%%%%%%%%%%%%%%%%%%%%
%\ \ \ \\
\begin{center}
\begin{minipage}{5.0in}
%%%%%%%%%%%%%%%%%%%%%%%%%%%%%%%%%%%%%%%%%%%%%%%%%%%%%%%%%%%%%%%%%%%%%%%%%%
\ \ \ \\
%\\
\begin{abstract}
\noindent Adomian decomposition method is used for solving the
seventh order boundary value problems. The approximate solutions
of the problems are calculated in the form of a rapid convergent
series and not at grid points. Two numerical examples have been
considered to illustrate the efficiency and implementation of the
method.
\end{abstract}
\end{minipage}
\end{center}
{\bf Keywords:} {\small
Adomian decomposition method; Seventh order boundary value problems; Linear and nonlinear problems; Series solution.} \\
%\end{keywords}
%%%%%%%%%%%%%%%%%%%%%%%%%%%%%%%%%%%%%%%%%%%%%%%%%%%%%%%%%%%%%%%%%%%%%%%%%%
%%%%%%%%%%%%%%%%%%%% This file is to create final form of septic %%%%%%%%%%
%\newpage
\newsection{Introduction}The boundary value problems of ordinary differential equations
play an important role in many fields. The theory of seventh order
boundary value problems is seldom in the numerical analysis
literature. These problems
 generally arise in modelling induction motors with two rotor
 circuits. The induction motor behavior is represented by a fifth
 order differential equation model. This model contains two stator state
 variables, two rotor state variables and one shaft speed.
Normally, two more variables must be added to account for the
effects of a second rotor circuit representing deep bars, a
starting cage or rotor distributed parameters. To avoid the
computational burden of additional state variables when additional
rotor circuits are required, model is often limited to the fifth
order and rotor impedance is algebraically altered as function of
rotor speed. This is done under the assumption that the frequency
of rotor currents depends on rotor speed. This approach is
efficient for the steady state response with sinusoidal voltage,
but it does not hold up during the transient conditions, when
rotor frequency is not a single value. So, the behavior of such
models show up in the seventh order boundary value problems
~\cite{I1}.\\Presently, the literature on the numerical solutions
of seventh order boundary value problems and associated eigen
value problems is seldom. Siddiqi and Ghazala ~\cite{SG5, SG6,
SG2, SG7, SG10, SG1, SG12} presented the solutions of fifth,
sixth, eighth, tenth and twelfth order boundary value problems
using polynomial and non-polynomial spline techniques. Siddiqi and
Twizell ~\cite{EH1996, EH1, EH2} presented the solution of eighth,
tenth and twelfth order boundary value problems using eighth,
tenth and twelfth degree splines respectively. \\
The Adomian decomposition method was introduced and developed by
George Adomian ~\cite{A1986, A1988, A1994, A1993}. This method has
been applied to a wide class of linear and nonlinear ordinary
differential equations, partial differential equations, integral
equations and integro-differential equations ~\cite{D2008,
Int1987, Int1998, Int2001}. Wazwaz~\cite{Int2001, In2001} provided
the solution of  fifth order and sixth order boundary value
problems by the modified decomposition method. Convergence of
Adomian decomposition method was studied by Cherruault $et~ al.$
~\cite{Int1992, Int1990, Int1993}. In this paper, the Adomian
decomposition method~\cite{A1986, A1988, A1994, A1993} is applied
on the seventh order boundary value problem. An algorithm for
approximate solution by Adomian method is developed for such
problems. This method
provides solutions in terms of a rapidly convergent series.  %This
\\
\newsection{Adomian Decomposition Method}
Consider the operator form of differential equation
\begin{eqnarray}\label{e2.1}
 Lu+ Ru + Nu= f,
\end{eqnarray}
 where $L$ is the highest order derivative operator with the assumption that $L^{-1}$ exists, $R$ is a linear derivative
operator of  order less than $L$, $Nu$ is the nonlinear term and
$f$ is the source term. Applying the inverse operator $L^{-1}$ to
both sides of Eq. (\ref{e2.1}) and using the given conditions
yields the following
\begin{eqnarray}\label{e2.2}
u= g - L^{-1}(Ru) - L^{-1} (Nu),
\end{eqnarray}
where the function $g$ represents the term obtained after
integrating the source term $f$ and using the given conditions .
The solution $u(x)$ can thus, be defined as
\begin{equation}\label{e2.3}
 u(x) = \sum^{\infty}_{n=0}u_{n}(x),
\end{equation}
where the components $u_{0}, u_{1}, u_{2}, ...$, are determined by
the following recursive relation
\begin{eqnarray}\label{e2.4}
u_{0}(x) = g,         ~~~~~~~~~~~~~~~~~~~~~~~~~~~~~~~~~~~~~~~~~~~ \nonumber  \\
 u_{k+1}(x) =  - L^{-1} R(u_{k}) - L^{-1}
N(u_{k}),~~k\geq0.
\end{eqnarray}
The nonlinear term $Nu$ is expressed in terms of an infinite
series of the Adomian polynomials given by
\begin{equation}\label{e2.5}
 Nu = \sum^{\infty}_{n=0}A_{n},
\end{equation}
where $A_{n}$ are Adomian polynomials  ~\cite{A1994} defined by
\begin{eqnarray}\label{e2.6}
A_{n}=\frac{1}{n!}\frac{d^{n}}{dp^{n}}\left[N\left(\sum^{\infty}_{k=0}p^{k}u_{k}
\right)\right]_{p=0}, ~~~~k=0, 1, 2, ...
\end{eqnarray}
\newsection{Analysis of the Method}
Consider the following seventh order boundary value problem
\begin{eqnarray}\label{e3.1}
 u^{(7)}(x)=\phi(x) u+ \psi(x) + f(x, u)  , ~~ 0 \leq x \leq b,
\end{eqnarray}
with boundary conditions
\begin{eqnarray}\label{e3.2}
\left.\begin{array}{rrr}
 u^{(i)}(0)&=& \alpha_{i},~~i=0,1,2,3,~~~~~~~~~~~~~~~~~~~~~~~~~~~~\\
 u^{(j)}(b)&=&\beta_{j},~~j=0,1,2.~~~~~~~~~~~~~~~~~~~~~~~~~~~~~~~
\end{array} \right\}\
 \end{eqnarray}
 where $\alpha_{i},~ i=0,1,2,3$ and $ \beta_{j},~
j=0,1,2$ are finite real constants and the function $f$ is
continuous on $[0,b]$.\\
 The problem (\ref{e3.1}), in operator form, can be rewritten as
\begin{eqnarray}\label{e3.3}
Lu =\phi(x)u+\psi(x) + f(x, u),
\end{eqnarray}
where the derivative operator $L$ is given by
\begin{eqnarray}\label{e3.4}
L=\frac{d^{7}}{dx^{7}} .
\end{eqnarray}
The inverse operator $L^{-1}$ is therefore defined as
\begin{eqnarray}\label{e3.5}
L^{-1}(.)=\int_{0}^{x}\int_{0}^{x}\int_{0}^{x}\int_{0}^{x}\int_{0}^{x}\int_{0}^{x}\int_{0}^{x}(.)~dx~
dx~dx~dx~dx~dx~dx.
\end{eqnarray}
Operating with $L^{-1}$ on Eq. (\ref{e3.3}) gives
\begin{eqnarray}\label{e3.6}
u(x)&=&\alpha_{0} +~ \alpha_{1}x + \frac{1}{2!}\alpha_{2}x^{2} +
\frac{1}{3!}\alpha_{3}x^{3} +\frac{1}{4!}\alpha_{4}x^{4}
+\frac{1}{5!}\alpha_{5}x^{5} +\frac{1}{6!}\alpha_{6}x^{6}
\nonumber
\\&&+ L^{-1}(\phi(x)u)+ L^{-1}(\psi(x)) + L^{-1}(f(x, u)),
\end{eqnarray}
using the boundary conditions (\ref{e3.2}) at  $x=0$ yields
\begin{eqnarray}\label{e3.7}
u(x)&=&\alpha_{0} + ~\alpha_{1}x + \frac{1}{2!}\alpha_{2}x^{2} +
\frac{1}{3!}\alpha_{3}x^{3} +\frac{1}{4!}A x^{4}
+\frac{1}{5!}Bx^{5} +\frac{1}{6!}Cx^{6} + \nonumber
\\&&+ L^{-1}(\phi(x)u)+ L^{-1}(\psi(x)) + L^{-1}(f(x, u)),
\end{eqnarray} where the
constants
\begin{eqnarray}\label{e3.8} A=u^{(4)}(0), ~~B=u^{(5)}(0),
~~C=u^{(6)}(0),
\end{eqnarray}
will be determined using boundary conditions at $x=b$.\\The
Adomian method determines the solution $u(x)$ in terms of the
following decomposition series
\begin{equation}\label{e3.9}
 u(x) = \sum^{\infty}_{n=0}u_{n}(x).
\end{equation}
Substituting Eq. (\ref{e3.9}) into Eq. (\ref{e3.7}) leads to
\begin{eqnarray}\label{e3.10}
\sum^{\infty}_{n=0}u_{n}(x)&=&\alpha_{0} + ~\alpha_{1}x +
\frac{1}{2!}\alpha_{2}x^{2} + \frac{1}{3!}\alpha_{3}x^{3}
+\frac{1}{4!}A x^{4} +\frac{1}{5!}Bx^{5} +\frac{1}{6!}Cx^{6} +
\nonumber
\\&&+ L^{-1}(\phi(x) \sum^{\infty}_{n=0}u_{n}(x))+ L^{-1}(\psi(x)) + L^{-1}(f(x,
u)).
\end{eqnarray}
To determine the components $u_{n}(x)$, $n\geq0$,  the following
recurrence relation will be used
\begin{eqnarray}\label{e3.11}
 u_{0}(x)&=&\alpha_{0} + ~\alpha_{1}x +
\frac{1}{2!}\alpha_{2}x^{2} + \frac{1}{3!}\alpha_{3}x^{3}
+\frac{1}{4!}A x^{4} +\frac{1}{5!}Bx^{5} +\frac{1}{6!}Cx^{6} +
 L^{-1}(\psi(x)), \nonumber\\
 u_{k+1}(x)&=&L^{-1}(\phi(x) \sum^{\infty}_{k=0}u_{k}(x))+  L^{-1}(\sum^{\infty}_{k=0}A_{k}), ~~~k\geq0,
 \end{eqnarray}
 where the nonlinear function $f (x, u)$ can be decomposed into an infinite series of Adomian polynomials given by
 \begin{eqnarray}\label{e3.12}
f(x, u)=\sum^{\infty}_{k=0}A_{k},
 \end{eqnarray}
 where $A_{k}$ are Adomian polynomials~\cite{A1994} defined by
\begin{eqnarray}\label{e3.13}
A_{k}=\frac{1}{k!}\frac{d^{k}}{dp^{k}}\left[N\left(\sum^{\infty}_{j=0}p^{j}u_{j}
\right)\right]_{p=0}, ~~~~n=0, 1, 2, ...
\end{eqnarray}
Keeping in view Eq. (\ref{e3.11}) - Eq. (\ref{e3.13}), the
components $u_{0}, u_{1}, u_{2}, ...$ are calculable. By the
Adomian
 method the solution can be constructed as
\begin{eqnarray*}
u=\lim_{n\rightarrow\infty}\varphi_{n},
\end{eqnarray*}
where the $n$-term approximant is defined by
\begin{eqnarray*}
\varphi_{n}=\sum^{n-1}_{i=0}u_{i}.
\end{eqnarray*}
Applying the boundary conditions at $x=b$ to the approximant
$\varphi_{n}$. The resulting algebraic system in $A$, $B$ and $C$
can be solved to produce approximations to constants $A$, $B$ and
$C$. Finally, the approximated solution of the seventh order
boundary value
problem follows.\\
 To implement the method, four
numerical examples are considered in the following section.
\newsection{Numerical Examples}
\textbf{Example 4.1}~~Consider the linear seventh order boundary
value problem
\begin{eqnarray}\label{e4.1}
 u^{(7)}(x)&=& x u(x)+e^x (x^2 - 2 x - 6),~~~~~~ 0 \leq  x \leq1,
 \end{eqnarray}
 subject to the boundary conditions
 \begin{eqnarray}\label{e4.2}
 \left.\begin{array}{lll}
u(0)&=&~~1,~~~~~u(1)=0,~~~\\
u^{(1)}(0)&=&~~0,~~u^{(1)}(1)=-e,~\\
u^{(2)}(0)&=&-1,~~u^{(2)}(1)=-2e,\\
u^{(3)}(0)&=&-2.~\\
\end{array} \right\}\
\end{eqnarray}
The exact solution of the problem (4.1)  is
 $$u(x)=(1-x)e^x.$$
 The problem (\ref{e4.1}), in operator form, can be rewritten as
\begin{eqnarray}\label{e4.3}
Lu &=&x u(x)+e^x(x^2 - 2 x - 6),
\end{eqnarray}
Operating with $L^{-1}$ on Eq. (\ref{e4.3}) and using the boundary
conditions (\ref{e4.2}) at  $x=0$ gives
\begin{eqnarray}\label{e4.4}
u(x)&=&-63 - 64x - \frac{35}{2!}x^{2} - 4x^{3}
+(-\frac{1}{2}+\frac{A}{24} )x^{4} +(-\frac{1}{30}+\frac{B}{120}
)Bx^{5} +\frac{1}{360}(2+C)x^{6} \nonumber \\ &&+e^x (-8 +x)^2+
L^{-1}(xu(x)),
\end{eqnarray}
where the constants
\begin{eqnarray}\label{e4.5}
A=u^{(4)}(0), ~B=u^{(5)}(0), ~C=u^{(6)}(0), ~
\end{eqnarray}
are to be determined. Substituting the decomposition series
(\ref{e2.3}) for $u(x)$ in Eq. (\ref{e4.4}) yields
\begin{eqnarray} \label{e4.6}
\sum^{\infty}_{n=0}u_{n}(x)&=&-63 - 64x - \frac{35}{2!}x^{2} -
4x^{3} +(-\frac{1}{2}+\frac{A}{24} )x^{4}
+(-\frac{1}{30}+\frac{B}{120} )Bx^{5} +\frac{1}{360}(2+C)x^{6}
\nonumber \\ &&+ e^x (-8 +x)^2+
L^{-1}(x\sum^{\infty}_{n=0}u_{n}(x)).
\end{eqnarray}
Using the recurrence algorithm (\ref{e3.11}) the following
relations are obtained
\begin{eqnarray}\label{e4.7}
 u_{0}(x)&=&-63 - 64x - \frac{35}{2!}x^{2} - 4x^{3}
+(-\frac{1}{2}+\frac{A}{24} )x^{4} +(-\frac{1}{30}+\frac{B}{120}
)Bx^{5} +\frac{1}{360}(2+C)x^{6} \nonumber \\
&&+ e^x (-8 +x)^2,\nonumber \\
 u_{k+1}(x)&=& L^{-1}(xu_{k}(x)), ~~~k\geq0.
 \end{eqnarray}
 From Eq. (\ref{e4.7}),
\begin{eqnarray}\label{e4.8}
 u_{0}(x)&=&-63 - 64x - \frac{35}{2!}x^{2} - 4x^{3}
+(-\frac{1}{2}+\frac{A}{24} )x^{4} +(-\frac{1}{30}+\frac{B}{120}
)Bx^{5} +\frac{1}{360}(2+C)x^{6} \nonumber \\
 && + e^x (-8 +x)^2,\nonumber \\
 u_{1}(x)&=& L^{-1}(xu_{0}(x)), \nonumber \\
 &=&1848 + 1392x + 505x^2 + 116x^3 + \frac{37}{2}x^4 +
 \frac{31}{15}x^5 +  \frac{17}{120}x^6 - \frac{x^8}{640} -
 \frac{x^9}{3780}-  \frac{x^{10}}{34560}  \nonumber \\&&-  \frac{x^{11}}{415800}
 + \frac{(-12 +A)x^{12}}{95800320}+ \frac{(-4 +B)x^{13}}{1037836800}+ \frac{(2 +C)x^{14}}{12454041600}
 + e^x(-1848 + 456x \nonumber \\&&- 37x^2 + x^3),\nonumber \\
 &\vdots&
 \end{eqnarray}
The series solution of $u(x)$ is found to be as an approximation
with  four components
 \begin{eqnarray}\label{e4.9}
 u(x)&=&u_{0}(x) + u_{1}(x)+u_{2}(x) + u_{3}(x).
  \end{eqnarray}
The unknown constants $A$, $B$ and $C$ can be obtained by imposing
the boundary conditions at $x = 1$ on the four-term approximant
$\varphi_{4}$ defined by Eq. (\ref{e4.9}) , as
$$A =-3.0000001004600083 ,~~~ B =-3.9999991664171546 ,~~~ C =-5.000002140841064
.$$ Finally, the series solution can be written as
\begin{eqnarray}\label{e4.10}
 u(x) &=&1 - \frac{x^{2}}{2} - \frac{x^{3}}{3} -0.125 x^4
 -0.0333333 x^{5} -0.00694445x^{6}-\frac{ x^{7}}{840}-\frac{ x^{8}}{5760}-\frac{ x^{9}}{45360}\nonumber \\
&& -\frac{ x^{10}}{403200}-\frac{ x^{11}}{3991680}
-(2.29644\times10^{-8})x^{12} - (1.60591 \times 10^{-10})x^{13} \nonumber \\
 &&+O(x^{14}).
 \end{eqnarray}
 The comparison of the exact solution
with the series solution of the problem (4.1) is given in Table 1,
which shows that the method is quite efficient. Figure 1 also
endorses the efficiency of the
method.\\\\
\textbf{Example 4.2}~~Consider the following seventh order nonlinear
boundary value problem
\begin{eqnarray}\label{e4.11}
 u^{(7)}(x) &=&-e^{x}u^{2}(x),~~ 0 <
x<1,
\end{eqnarray}
subject to the boundary conditions
\begin{eqnarray}\label{e4.12}
 \left.\begin{array}{lll}
u(0)&=&~~1,~~~~~~u(1)=e^{-1},~~~\\
u^{(1)}(0)&=&-1,~~~u^{(1)}(1)=-e^{-1},~\\
u^{(2)}(0)&=&~~1,~~~u^{(2)}(1)=e^{-1},~~~\\
u^{(3)}(0)&=&-1.~\\
\end{array} \right\}\
\end{eqnarray}
The exact solution of the problem (4.2) is
 $$u(x)=e^{-x}.$$
The problem (\ref{e4.11}), in operator form, can be rewritten as
\begin{eqnarray}\label{e4.13}
Lu &=&e^{x}u^{2}(x),
\end{eqnarray}
Operating with $L^{-1}$ on Eq. (\ref{e4.13}) and using the
boundary conditions (\ref{e4.12}) at  $x=0$ gives
\begin{eqnarray}\label{e4.14}
u(x)=1 - x + \frac{1}{2!}x^{2} - \frac{1}{3!}x^{3} +\frac{1}{4!}A
x^{4} +\frac{1}{5!}Bx^{5} +\frac{1}{6!}Cx^{6} +
L^{-1}(-e^{x}u^{2}(x)),
\end{eqnarray}
where the constants
\begin{eqnarray}\label{e4.15}
A=u^{(4)}(0), ~B=u^{(5)}(0), ~C=u^{(6)}(0), ~
\end{eqnarray}
are to be determined. Substituting the decomposition series
(\ref{e2.3}) for $u(x)$ and the series of polynomials (\ref{e2.5})
for $u^2(x)$ in Eq. (\ref{e4.14}) yields
\begin{eqnarray}\label{e4.16}
\sum^{\infty}_{n=0}u_{n}(x)=&1& - x + \frac{1}{2!}x^{2} -
\frac{1}{3!}x^{3} +\frac{1}{4!}A x^{4} +\frac{1}{5!}Bx^{5}
+\frac{1}{6!}Cx^{6} \nonumber \\ &+&
L^{-1}(-e^{x}\sum^{\infty}_{n=0}A_{n}),
\end{eqnarray}
where $A_{n}$ are the Adomian polynomials for the nonlinear term
$Nu= u^{2}(x)$, thus given by
\begin{eqnarray}\label{e4.17}
A_{0}&=&N(u_{0})= u_{0}^{2}(x),\nonumber \\
A_{1}&=&u_{0}(x) N^{'}(u_{0})= 2u_{0}(x)u_{1}(x),\nonumber \\
A_{2}&=&u_{2}(x) N^{'}(u_{0})+\frac{u^{2}_{1}(x)}{2!}N^{''}(u_{0})= 2u_{0}(x)u_{2}(x)+u^{2}_{1}(x),\nonumber \\
&\vdots&
\end{eqnarray}
Using the recurrence algorithm (\ref{e3.11}) yields
\begin{figure}
\centering
\begin{tabular}{cc}
\begin{minipage}{200pt}
\frame{\includegraphics[width=200pt]{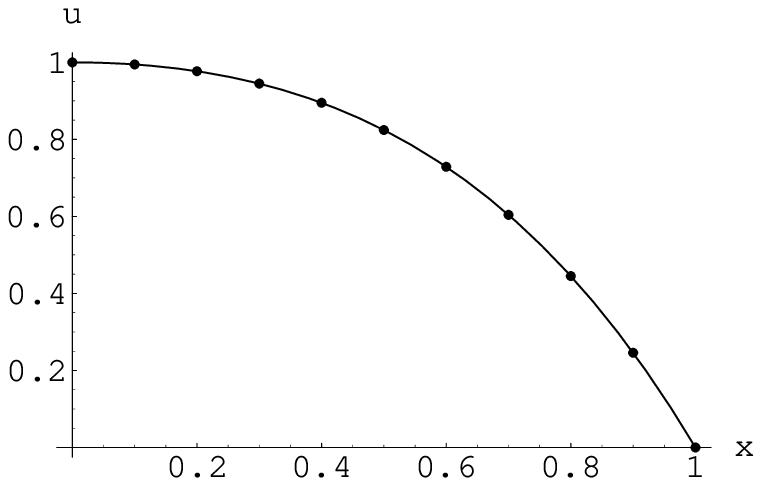}} \caption{Comparison
of the approximate solution with the exact solution for the problem
(4.1). Dotted line: approximate solution, solid line: the exact
solution.}
\end{minipage}
&
\begin{minipage}{200pt}
\frame{\includegraphics[width=200pt]{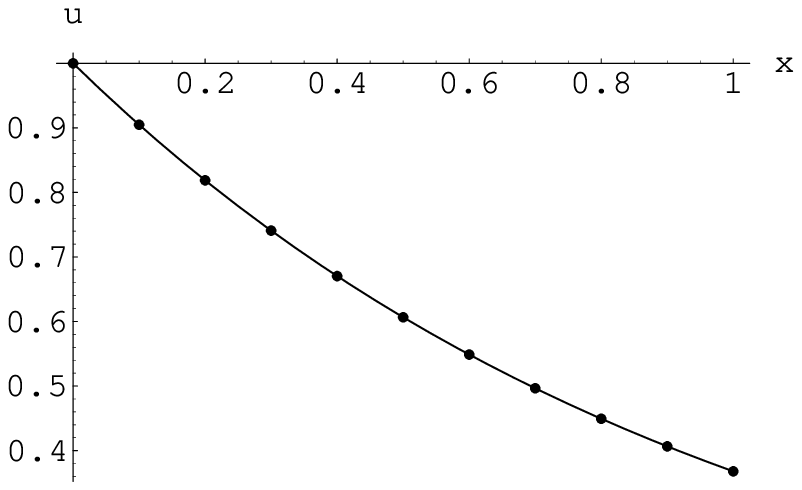}} \caption{Comparison
of the approximate solution with the exact solution for the problem
(4.2). Dotted line: approximate solution, solid line: the exact
solution.}
\end{minipage}
\end{tabular}
\end{figure}\\
\begin{eqnarray}\label{e4.18}
 u_{0}(x)&=&1 - x + \frac{1}{2!}x^{2} - \frac{1}{3!}x^{3} +\frac{1}{4!}A
x^{4} +\frac{1}{5!}Bx^{5} +\frac{1}{6!}Cx^{6}, \nonumber\\
 u_{k+1}(x)&=& L^{-1}(-e^{x}A_{k}), ~~~k\geq0,
 \end{eqnarray}
 From (\ref{e4.17}) and (\ref{e4.18}),
 \begin{eqnarray*}
 u_{0}(x)&=&1 - x + \frac{1}{2!}x^{2} - \frac{1}{3!}x^{3} +\frac{1}{4!}A
x^{4} +\frac{1}{5!}Bx^{5} +\frac{1}{6!}Cx^{6}, \\
 u_{1}(x)&=& L^{-1}(-e^{x}A_{0}),\\
&=&- \frac{x^{7}}{5040} + \frac{x^{8}}{40320} -
\frac{x^{9}}{362880} + \frac{x^{10}}{3628800} +
(\frac{1}{39916800}-\frac{A}{19958400})x^{11} \\ &&-
(\frac{1}{479001600}+\frac{B}{239500800})x^{12}+
(\frac{1}{6227020800}-\frac{C}{3113510400})x^{13} +O(x^{14}),
 \end{eqnarray*}
 The series solution of $u(x)$ is found to be as an approximation with  two components
 \begin{eqnarray}\label{e4.19}
 u(x)&=&u_{0}(x) + u_{1}(x).
  \end{eqnarray}
 The unknown constant $A$, $B$, and $C$ can be obtained by imposing the
boundary conditions at $x = 1$ on the two-term approximant
$\varphi_{2}$ defined by Eq. (\ref{e4.19}) , as
$$A =1.0000000197456873 ,~~~ B =-1.000000228700079 ,~~~ C =1.0000008112153262
.$$ Finally, the series solution can be written as
\begin{eqnarray} \label{e4.20}
 u(x) &=&1 - x + 0.5x^{2} - 0.166667x^{3} +0.0416667
x^{4} -0.00833334x^{5} + 0.00138889x^{6} \nonumber \\
&& - 0.000198413x^{7} + 0.0000248016x^{8}- (2.75573 \times
10^{-6})x^{9} + (2.75573 \times 10^{-7})x^{10} \nonumber \\
&& - (2.50521 \times 10^{-8})x^{11} +(2.08768 \times
10^{-9})x^{12} - (1.60591 \times 10^{-10})x^{13} \nonumber \\
 &&+O(x^{14}).
 \end{eqnarray}
  The comparison of the exact
solution with the solution of the problem (4.2) is given in Table 2,
which shows that the method is quite efficient. Figure 2 also
endorses the efficiency of the method.\\\begin{figure} \centering
\begin{tabular}{cc}
\begin{minipage}{200pt}
\frame{\includegraphics[width=200pt]{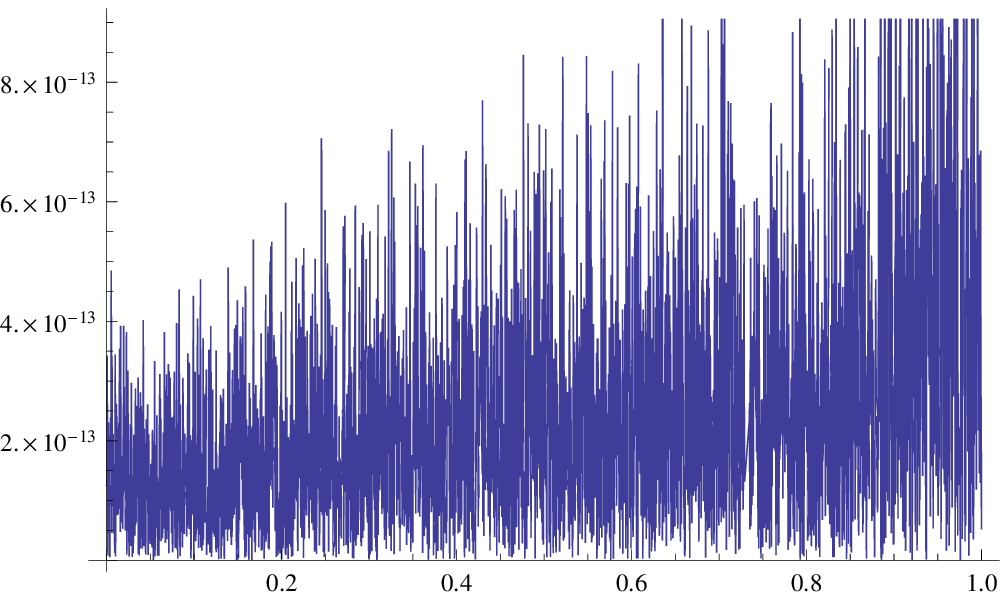}} \caption{Absolute
Error for the problem (4.3). }
\end{minipage}
&
\begin{minipage}{200pt}
\frame{\includegraphics[width=200pt]{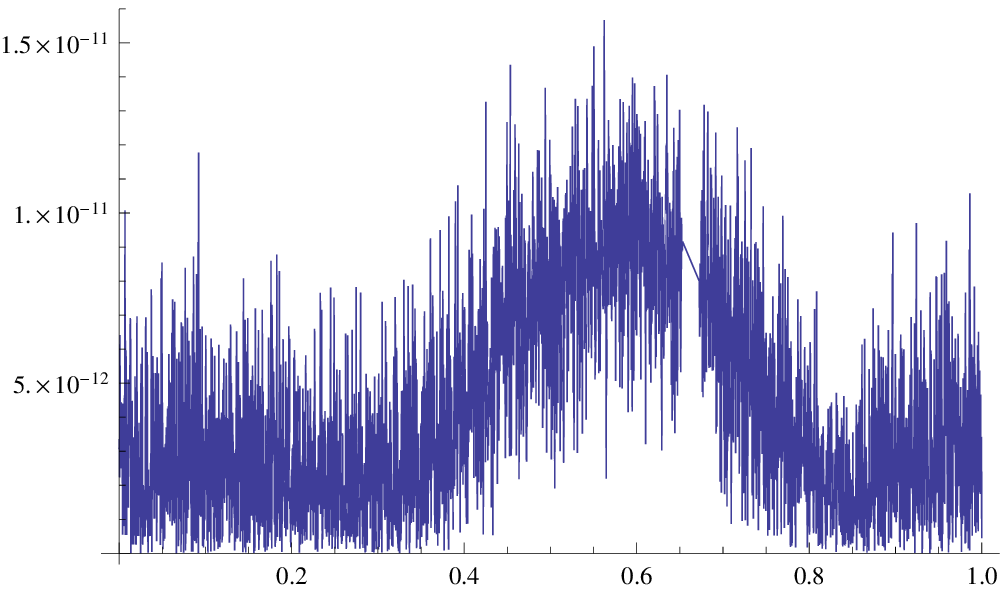}} \caption{Absolute
Error for the problem (4.4).}
\end{minipage}
\end{tabular}
\end{figure}\\
\textbf{Example 4.3}~~Consider the following seventh order linear
boundary value problem
\begin{eqnarray}\label{e4.21}
\left.\begin{array}{lll}
 u^{(7)}(x)&=&-u(x)-e^x (35+12x+2x^2),~~~~ 0 \leq  x \leq1,\\\quad
u(0)&=&~0,~~~~~~u(1)=0,~~~\\
u^{(1)}(0)&=&~1,~~~u^{(1)}(1)=-e,~\\
u^{(2)}(0)&=&~0,~~~u^{(2)}(1)=-4e,\\
u^{(3)}(0)&=&-3.~~\\
\end{array} \right\}\
\end{eqnarray}
The exact solution of the Example 3.1  is
 $u(x)=x(1-x)e^x $ \cite{Muza12}.\\
 Following the procedure of the Example 4.1, this problem is
 solved. It is observed that the errors in absolute values are better than
those of Siddiqi and Iftikhar \cite{Muza12} as shown in Table 3. In
Figure 3 absolute errors are plotted.
 \\\\
 \textbf{Example 4.4}~~The following seventh order nonlinear boundary
value problem is considered
\begin{eqnarray}\label{e4.22}
\left.\begin{array}{lll}
 u^{(7)}(x)&=&u(x) u'(x)+e^{-2x} (2+e^x(x-8)-3x+x^2),~~~~ 0 \leq  x \leq1,\\\quad
u(0)&=&1,~~~~~u(1)=0,~~~\\
u^{(1)}(0)&=&-2,~~u^{(1)}(1)=-1/e,~\\
u^{(2)}(0)&=&3,~~u^{(2)}(1)=2/e,\\
u^{(3)}(0)&=&-4.~~\\
\end{array} \right\}\
\end{eqnarray}
The exact solution of the Example 3.2  is
 $u(x)=(1-x)e^{-x}. $\\
  Following the procedure of the previous problem 4.2, this problem is
 solved. The comparison of the exact
solution with the solution of the problem 4.2 is given in Table 4,
which shows that the method is quite efficient. In Figure 4
absolute errors are plotted.\\
\textbf{Conclusion}~~In this paper, the Adomian decomposition method
has been applied to obtain the numerical solutions of linear and
nonlinear seventh order boundary value problems.  The numerical
results show that the method is quite efficient for solving high
order boundary value problems arising in various fields of
engineering and science.
\\
\begin{table}[bht]
\caption{Comparison of numerical results for the problem 4.1}
\centering
\begin{small}
\begin{tabular}{|c|c|c|c|}
\hline $x$& Exact solution  & Approximate series solution & Absolute Error \\
\hline$0.0$&1.0000&1.0000&0.0000\\
\hline$0.1$&0.0994&0.0994&4.3972E-10\\
\hline$0.2$&0.9771& 0.9771&4.9251E-10\\
\hline$0.3$&0.9449& 0.9449& 7.4067E-10\\
\hline$0.4$&0.8950& 0.8950&6.6537E-10\\
\hline$0.5$&0.8243& 0.8243&3.0059E-11\\
\hline$0.6$&0.7288& 0.7288& 4.3591E-10\\
\hline$0.7$&0.6041& 0.6041&3.6735E-10\\
\hline$0.8$&0.4451&0.4451&7.2753E-10\\
\hline$0.9$&0.2459&0.2459&7.0036E-10\\
\hline$1.0$&0.0000&2.2191E-10&2.2191E-10\\
\hline
\end{tabular}\end{small}
\end{table}
\begin{table}[bht]
\caption{Comparison of numerical results for problem 4.2} \centering
\begin{small}
\begin{tabular}{|c|c|c|c|}
\hline $x$& Exact solution  & Approximate series solution & Absolute Error \\
\hline$0.0$&1.0000&1.0000&0.0000\\
\hline$0.1$&0.9048&0.9048&1.5676E-9\\
\hline$0.2$&0.8187& 0.8187&1.6418E-9\\
\hline$0.3$&0.7408& 0.7408& 4.9680E-9\\
\hline$0.4$&0.6703& 0.6703&1.5514E-9\\
\hline$0.5$&0.6065& 0.6065&1.5274E-9\\
\hline$0.6$&0.5488&0.5488& 2.4958E-9\\
\hline$0.7$&0.4965& 0.4965&1.3993E-8\\
\hline$0.8$&0.4493& 0.4493&2.5593E-9\\
\hline$0.9$&0.4065&0.4065&5.4089E-9\\
\hline$1.0$&0.3678&0.3678&1.1034E-9\\
\hline
\end{tabular}\end{small}
\end{table}
\begin{table}[bht]
\caption{Comparison of numerical results for the problem 4.3}
\centering
\begin{small}
\begin{tabular}{|c|c|c|c|c|}
\hline $x$& Exact   & Approximate   &Absolute Error &Absolute Error \\
 & solution  & series solution &present method& Siddiqi and Iftikhar \cite{Muza12} \\
\hline$0.0$&0.0000&0.0000&0.0000&0.0000\\
\hline$0.1$&0.9946&0.9946&1.23082E-13&8.55607E-13\\
\hline$0.2$&0.1954&0.1954&3.7792E-13&9.94041E-12\\
\hline$0.3$&0.2835& 0.2835&2.37421E-13& 3.52244E-11\\
\hline$0.4$&0.3580& 0.3580&3.62099E-13&7.3224E-10\\
\hline$0.5$&0.4122& 0.4122&9.39249E-14&1.08769E-10\\
\hline$0.6$&0.4373& 0.4373&4.82947E-13& 1.29035E-10\\
\hline$0.7$&0.4229& 0.4229&1.09135E-13&1.51466E-10\\
\hline$0.8$&0.3561& 0.3561&1.64868E-14&2.717974E-10\\
\hline$0.9$&0.2214&0.2214&7.25975E-13&7.48179E-10\\
\hline$1.0$&0.0000&-4.54747E-13&4.54747E-13&2.1729E-09\\
\hline
\end{tabular}\end{small}
\end{table}
\begin{table}[bht]
\caption{Comparison of numerical results for the problem 4.4}
\centering
\begin{small}
\begin{tabular}{|c|c|c|c|}
\hline $x$& Exact solution  & Approximate series solution & Absolute Error \\
\hline$0.0$&1.0000&1.0000&1.67932E-12\\
\hline$0.1$&0.814354&1.1051&2.96696E-12\\
\hline$0.2$&0.654985& 0.654985&1.26055E-12\\
\hline$0.3$&0.518573& 0.518573& 2.10898E-12\\
\hline$0.4$&0.402192& 0.402192&6.68926E-12\\
\hline$0.5$&0.303265& 0.303265&7.21923E-12\\
\hline$0.6$&0.219525& 0.219525& 9.75339E-12\\
\hline$0.7$&0.148976& 0.148976&2.19552E-12\\
\hline$0.8$&0.0898658& 0.0898658&4.24917E-12\\
\hline$0.9$&0.040657&0.040657&2.27311E-13\\
\hline$1.0$&0.0000&4.42298E-12&4.42298E-12\\
\hline
\end{tabular}\end{small}
\end{table}
\bibliographystyle{amsplain}
\bibliography{xbib}
\end{document}